\documentclass[12pt]{article} 

\usepackage[utf8]{inputenc} 
\usepackage{bbm}
\usepackage{relsize}
\usepackage{authblk}
\usepackage{geometry} 
 \usepackage{amsthm}
\usepackage{
	amsfonts,
	amsmath,
	appendix,
	color,
	amssymb,
	graphicx,
	bm,
	stmaryrd,
	mathabx,
	amssymb,
	mathdots,
	esint
}
\usepackage{amscd,mathrsfs,pifont}
\usepackage{eucal,color,empheq}

\newtheorem{thm}{Theorem}
\newtheorem{lemma}{Lemma}

\newcommand{\bea}{\begin{eqnarray}}
\newcommand{\eea}{\end{eqnarray}}
\newcommand{\R}{\mathbb{R}}
\newcommand{\simb}{\mathlarger{\mathlarger{\sim}}}
\newcommand{\tr}{\mbox{Tr}}
\DeclareMathOperator{\erfc}{erfc}
\DeclareMathOperator{\erfcx}{erfcx}
\def\sumd{\displaystyle\sum}

\newcommand{\de}[1]{d #1}

\newcommand{\der}[1]{\frac{d}{d #1}}

\newcommand{\Obig}[1]{\underline{O}\left(#1\right)}

\newcommand{\expp}[2]{e^{\left(+\right)}_{#1}\left(#2\right)}
\newcommand{\expm}[2]{e^{\left(-\right)}_{#1}\left(#2\right)}
\title{What is the probability that a large random matrix has no
real eigenvalues?}
\author[1]{Eugene Kanzieper\thanks{eugene.kanzieper@hit.ac.il}}
\author[2]{Mihail Poplavskyi\thanks{m.poplavskyi@warwick.ac.uk}}
\author[3]{Carsten Timm\thanks{carsten.timm@tu-dresden.de}}
\author[2]{\newline Roger Tribe\thanks{r.p.tribe@warwick.ac.uk}}
\author[2]{Oleg Zaboronski\thanks{olegz@maths.warwick.ac.uk}}
\affil[1]{\small{Department of Applied Mathematics, Holon Institute of Technology, Holon 5810201, Israel and Department of Physics of Complex Systems, Weizmann Institute of Science, Rehovot 7610001, Israel}}
\affil[2]{\small{Department of Mathematics, University of Warwick, Coventry CV4~7AL, UK}}
\affil[3]{\small{Institute of Theoretical
Physics, Technische Universit\"at Dresden, 01062 Dresden, Germany}}

\date{\today}

\begin{document}
\maketitle
\begin{abstract}
We study the large-$n$ limit of the probability $p_{2n,2k}$ that a random $2n\times 2n$
matrix sampled from the real Ginibre ensemble has $2k$ real
eigenvalues. 
We prove that, $$\lim_{n\rightarrow \infty}\frac {1}{\sqrt{2n}} \log p_{2n,2k}=\lim_{n\rightarrow \infty}\frac {1}{\sqrt{2n}} \log p_{2n,0}=
-\frac{1}{\sqrt{2\pi}}\zeta\left(\frac{3}{2}\right),$$ where $\zeta$ is the Riemann zeta-function. 
Moreover, for any sequence of non-negative integers $(k_n)_{n\geq 1}$,  $$\lim_{n\rightarrow \infty}\frac {1}{\sqrt{2n}} \log p_{2n,2k_n}=-\frac{1}{\sqrt{2\pi}}\zeta\left(\frac{3}{2}\right),$$
provided $\lim_{n\rightarrow \infty} \left(n^{-1/2}\log(n)\right) k_{n}=0$.
\end{abstract}
\section{Introduction and the main result.}
Our paper is dedicated to the study of  the probability $p_{2n,2k}$ that a real
$2n\times 2n$ random matrix with independent normal entries (the so called `real Ginibre matrix') has
$2k$ real eigenvalues.
It has been known since \cite{edelman} that a typical large $N\times N$ Ginibre
matrix  has $O(\sqrt{N})$ real
eigenvalues. What is the probability of rare events consisting
of such a matrix having either anomalously many or few real eigenvalues? 

The former question has been addressed by many authors. Building on
the original work \cite{ginibre} by Ginibre, Edelman used the real Schur decomposition to prove that
\bea
p_{N,N}=\left(\frac{1}{2}\right)^{\frac{N(N-1)}{4}} \nonumber
\eea
see \cite{edelman2}.
In \cite{eugene1} Akemann and Kanzieper  employed the method of skew-orthogonal polynomials
to determine the probability that
all but two eigenvalues of a real 
Ginibre matrix are real.  In the large-$N$ limit, their result reads
\bea
p_{N,N-2}=e^{- \frac{\log(2)}{4} N^2 + \frac{\log(3\sqrt{2})}{2} N+\underline{o}(N)},
\eea
where $\lim_{N\rightarrow \infty} \underline{o}(N)/N=0$.
These answers were generalized in  a very recent paper \cite{ldon} 
where the large deviations principle of 
\cite{benarous} was extended to prove that the probability that a real Ginibre matrix has $\alpha N$ (where $0<\alpha<1$)
real eigenvalues is
$
p_{N,\alpha N} \stackrel{N \rightarrow \infty}{\simb} e^{-N^2 I_\alpha},
$
where the symbol ''$\sim$'' denotes the logarithmic asymptotic equivalence and the constant $I_\alpha$
is characterised as the minimal value of an explicitly given rate functional, see Proposition 2 and formula (4) of \cite{ldon}. 

In the present paper we answer the question about the probability that
a real Ginibre matrix has very few real eigenvalues:
 
\begin{thm}\label{thm1}
Let $G_{2n}$ be a random $2n\times 2n$ real matrix with independent $N(0,1)$
matrix elements. Let $p_{2n,2k}$ be the probability
that $G_{2n}$ has $2k$ real eigenvalues.
Then for any fixed $k=0,1,2,3,\ldots$,
\begin{eqnarray}\label{result}
\lim_{n\rightarrow\infty} \frac{1}{\sqrt{2n}} \log p_{2n,2k} =-\frac{1}{\sqrt{2\pi}}\zeta\left(\frac{3}{2}\right), 
\end{eqnarray}
where $\zeta$ is the Riemann zeta-function. 
Moreover, 
\begin{eqnarray}\label{result2}
\lim_{n\rightarrow\infty} \frac{1}{\sqrt{2n}} \log p_{2n,2k_{n}} =-\frac{1}{\sqrt{2\pi}}\zeta\left(\frac{3}{2}\right),
\end{eqnarray}
where $(k_n)_{n\geq 1}$ is a sequence of non-negative integers such that
\bea
\lim_{n\rightarrow \infty} \left(n^{-1/2}\log(n)\right)k_n=0. \nonumber
\eea
\end{thm}

In particular, the probability that a large $2n\times 2n$ Ginibre matrix
has no real eigenvalue behaves as
\begin{eqnarray}
p_{2n,0}\stackrel{n\rightarrow \infty}{\simb} e^{-\sqrt{\frac{n}{\pi}}\zeta(3/2)+\underline{o}(\sqrt{n})}. \nonumber
\end{eqnarray}

Notice that the answer (\ref{result}) is qualitatively different from
the results for the probability of having $O(n)$ real eigenvalues quoted above:
the `cost' of having $O(n)$ real eigenvalues normalised by the total number
of `anomalous' eigenvalues increases linearly with $n$, whereas the `cost'
of removing all real eigenvalues from the real axis is constant per eigenvalue.

It is also worth noting that our result `almost' extends to the  typical region $k\sim n^{1/2}$. 
(For example, we can choose $k_n=[\sqrt{n}/\log^2{n}]$ in (\ref{result2}).)
It would be interesting to see if (\ref{result2}) survives for $k_n=[c\sqrt{n}]$ where $c\ll1$.

The statement of the theorem can be {\it guessed} using existing
results:
in the limit $N\rightarrow \infty$, the unscaled law of real eigenvalues for the real Ginibre $N\times N$
ensemble converges. The limit coincides with the $t=1$ law for the 
$A+A\rightarrow \emptyset$ interacting particle
system on $\mathbb{R}$ \cite{roger_oleg}. The probability that an interval
of length $s$ has no particles for $A+A\rightarrow \emptyset$
has been calculated formally by Derrida and Zeitak \cite{derrida}. These
two facts
allowed Forrester \cite{forrester} to conclude that the large-$N$ limit of the probability that
there are no real eigenvalues in the interval $(a,a+s)$ should be given by
\begin{eqnarray} 
Prob[G_\infty \mbox{ has no eigenvalues in } (a,a+s)]\stackrel{s\rightarrow \infty}{\simb} e^{-\frac{1}{2\sqrt{2\pi}}\zeta(3/2)s}.
\label{eq:forrester}
\end{eqnarray}
Let us stress that equation (\ref{eq:forrester}) is valid for $N=\infty$ only.
However, we know from the work of Borodin and Sinclair \cite{borodin_sinclair} and Forrester and Nagao \cite{forrester_nagao} that the law of real eigenvalues
for the real Ginibre ensemble is a Pfaffian point process for all values
of $N\leq \infty$. Convergence of the finite-$N$ kernel to
the $N=\infty$ kernel
is exponentially fast within the spectral radius. The spectral radius is $R_N=\sqrt{N}+O(1)$ \cite{edelman}. 
We also know that the boundary effects
for a large but finite matrix size $N$ are only felt in the boundary
layer of the width of order $1$ near the edge.  
Therefore, the simplest finite-$N$ guess for
 $Prob[G_{N}$ has no real eigenvalues$]$ is
\begin{eqnarray*}
&&Prob[G_{N}~\mbox{has no real eigenvalues}]\\
&\approx& Prob[G_{N}~\mbox{has no real eigenvalues in }
(-R_N+L,R_N-L)]\\
&\approx& Prob[G_\infty~\mbox{has no real eigenvalues in } (-R_N,R_N)]
\end{eqnarray*}
Here $L\gg1$ is a large $N$-independent constant. 
The last probability in our heuristic chain of arguments
can be approximated using (\ref{eq:forrester}) with $s=2R_N$.
This suggests
\begin{eqnarray} 
Prob[G_{N} \mbox{ has no real eigenvalues}]\approx e^{-\frac{1}{\sqrt{2\pi}}\zeta(3/2)\sqrt{N}},
\nonumber 
\end{eqnarray}
which agrees with the statement of Theorem \ref{thm1}.

The value of the constant which defines the rate of decay of 
$p_{2n,0}$ in (\ref{result})
is $$\frac{1}{\sqrt{2\pi}}\zeta(3/2)\approx  1.0422,$$
which is consistent with its numerical estimate, see Figure \ref{fig.pn0}. 
The numerical analysis of the exact formula for $p_{2n,0}$ (see (\ref{lm1f2}) below) also shows that under
the assumption that the next-to-leading
term in the large-$N$ expansion of $p_{N,0}$ is constant, the resulting coefficient ($\approx 0.06268$)
is close to its exact counterpart from the large gap size expansion
of the Derrida-Zeitak formula ($\approx 0.0627$). 
At the moment we do not have a theory 
explaining this closedness.  
\begin{figure}
\centering
\includegraphics[trim=0cm 0cm 0cm 0.0cm, clip=true, width=0.75\textwidth]{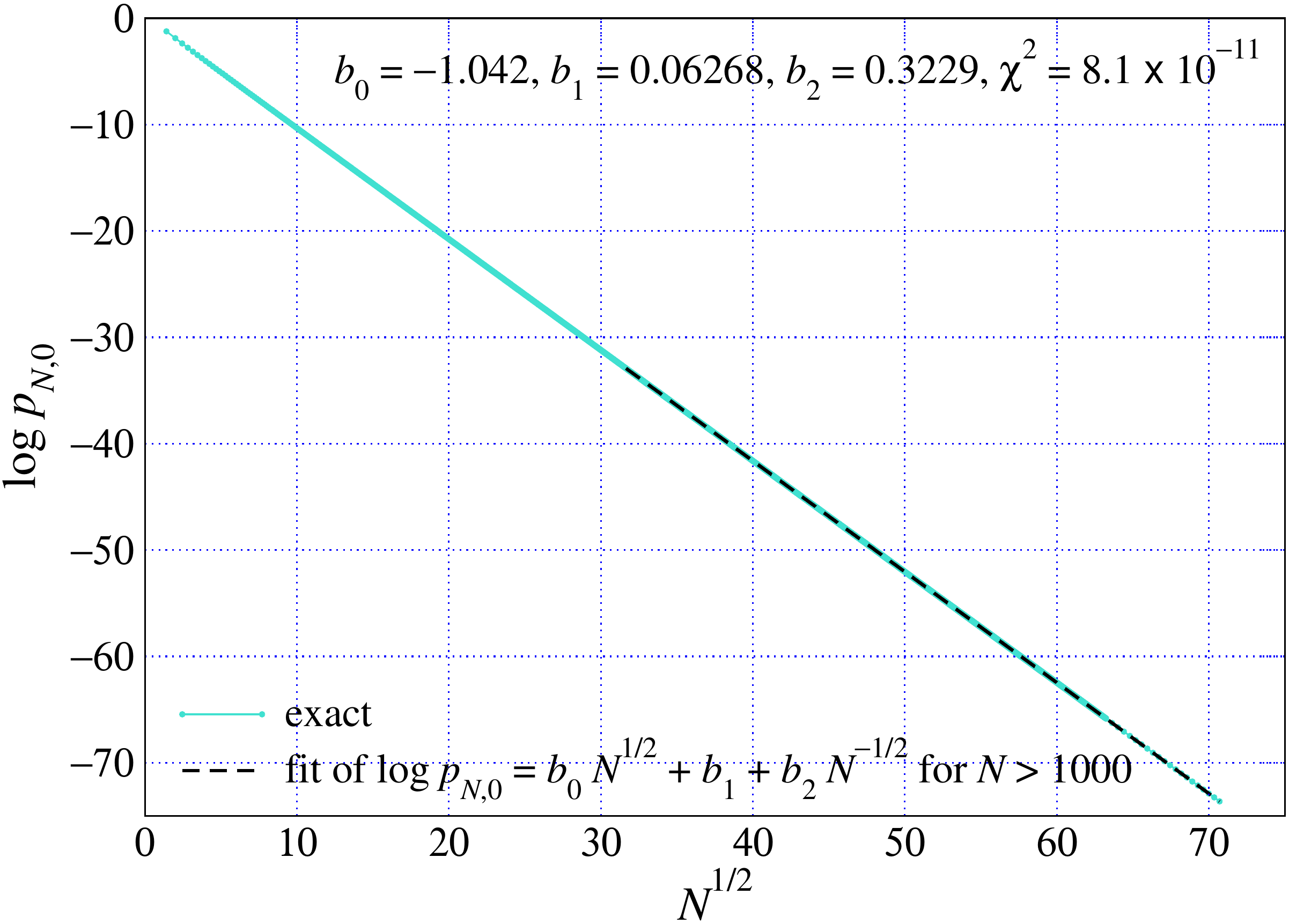}
\caption{The logarithm of the probability $p_{N,0}$ that an $N\times N$ 
matrix of even size sampled from the
real Ginibre ensemble does not have any real eigenvalues, as a function of
$\sqrt{N}$. The leading coefficient extracted using the best fit is
$-1.042$, the best fit for the next-to-leading constant is $0.06268$.
The `exact' curve is constructed using formula (\ref{lm1f2}) of Lemma \ref{lm1}
below. The form of the $b_2$-term in the fitting curve was chosen to minimize    the numerical goodness-of-fit $\chi^2$.
\label{fig.pn0}}
\end{figure}

Both the numerical simulations and the  heuristic argument given above
provide a strong hint in favour of Theorem \ref{thm1}.

There are several possible routes to the proof of the theorem. For example, one can try
to use Forrester's observation, coupled with
the knowledge of the rate of convergence of the 
Borodin-Sinclair-Forrester-Nagao kernel in the large-$N$
limit, to show that the errors in applying Derrida-Zeitak's formula
to gaps of  $N$-dependent sizes vanish as $N\rightarrow \infty$.
There is however a problem with this approach: in the case
we are interested in (annihilating Brownian motions or the $2$-state Potts model)
the infinite sums entering
 the gap
 formula converge only polynomially, see \cite{derrida}
for details.
Therefore, a careful 
justification would be required for the validity of the interchange of summation and taking
the large gap size limit. We feel that such a justification is best done
in the context of a general theory of `Fredholm Pfaffians'. In this paper, 
we will adopt {\em the spirit} of Derrida-Zeitak's calculation to construct rigorous asymptotics
of a very compact and easy to use exact determinantal
expression for the probability $p_{2n,2k}$ specific to the real Ginibre
ensemble. This determinantal expression 
can be derived building upon the results of \cite{KA-2005} and \cite{forrester_nagao}, see Lemma \ref{lm1}
below. We hope of course that our very
specialized proof will contribute to the general discussion of the theory 
of large deviations for Pfaffian point processes.

There is a drawback to our approach as well: even though we can now claim that (\ref{result}) is true, we still do not know how a large Ginibre matrix without real eigenvalues {\it looks}. For example, is there a unique optimal
configuration of complex eigenvalues for such matrices? 
What can be said about the overlaps between left and right eigenvectors of Ginibre matrices without real eigenvalues?
To answer these questions, one has to develop a large deviations principle along the lines of \cite{ldon} 
which will most likely use the picture of the `two-component' plasma consisting of one-dimensional and two-dimensional `gases' of eigenvalues discussed there.  

Our paper is organised as follows: a reader who is satisfied by our heuristic argument and the numerics can stop here. Those interested in the mathematical proof are advised
to read Section $2$ and consult Appendix \ref{app1} for the proofs of the technical
facts used in the proof of Theorem \ref{thm1}. Appendix \ref{app2} contains
remarks on the numerical evaluation of $p_{2n,0}$ for large values of $n$.

\section{The proof of Theorem \ref{thm1}.}

Our starting point is the following exact determinantal representation 
for the generating function for the probabilities $p_{2n,2k}$:
\begin{lemma}\label{lm1}
Let $n$ be a positive integer. Then
\begin{eqnarray}
\sum_{k=0}^n z^k p_{2n,2k}=\det_{j,k=1,n}\left[\delta_{j,k}+\frac{(z-1)}{\sqrt{2\pi}}\frac{\Gamma(j+k-3/2)}{\sqrt{\Gamma(2j-1)\Gamma(2k-1)}}\right].
\end{eqnarray}
In particular,
\begin{eqnarray} \label{lm1f2}
p_{2n,0}=\det_{j,k=1,n}\left[\delta_{j,k}-\frac{1}{\sqrt{2\pi}}\frac{\Gamma(j+k-3/2)}{\sqrt{\Gamma(2j-1)\Gamma(2k-1)}}\right].
\end{eqnarray}
\end{lemma}
We postpone to Appendix A the proofs of all lemmas used during the proof of the main theorem. \\

Notice that the expression (\ref{lm1f2}) coincides (as it should) with the $s\rightarrow \infty$
limit of the probability that a $2n\times 2n$ real Ginibre matrix has no
real eigenvalues in the interval $(-s,s)$ calculated by Forrester, see
formula (3.58) of \cite{forrester}.

We will prove Theorem \ref{thm1} in two steps: first, we will prove (\ref{result}) for $k=0$,
then we will show that $\lim_{n\rightarrow \infty} \frac{1}{\sqrt{2n}} \log p_{2n,2k_n}=
\lim_{n\rightarrow \infty} \frac{1}{\sqrt{2n}} \log p_{2n,0}$, where $(k_n)_{n \geq 1}$
is a sequence of integers which grows with $n$ slower than $n^{1/2}/\log(n)$.
\subsection{The calculation of $\lim_{n\rightarrow \infty} \frac{1}{\sqrt{2n}} \log p_{2n,0}$.}
Let $M_n$ be an $n\times n$ symmetric matrix entering the statement of Lemma \ref{lm1}: 
\bea \label{mat}
M_n(j,k)=\frac{1}{\sqrt{2\pi}}\frac{\Gamma(j+k-3/2)}{\sqrt{\Gamma(2j-1)\Gamma(2k-1)}},~1\leq j,k\leq n.
\eea
\begin{lemma}\label{lm2}
$M_n$ is a positive definite matrix. Moreover, there exists a positive constant
$\mu>0$ and a natural number $N$ such that for any $n>N$,
\bea \label{specrad}
\lambda_{max}(n)\leq 1-\frac{\mu}{n},
\eea 
where $\lambda_{max}(n)$ is the maximal eigenvalue of $M_n$.
\end{lemma}
Using Lemmas  \ref{lm1} and \ref{lm2} we represent $p_{2n,0}$ as follows:
\bea
\frac{1}{\sqrt{2n}}\log p_{2n,0} &=& \frac{1}{\sqrt{2n}} \tr \log(I-M_n) \nonumber \\
& = & -
\frac{1}{\sqrt{2n}}\sum_{m=1}^{K_n} \frac{1}{m} \tr M_n^m-\frac{1}{\sqrt{2n}}R_n(K_n), \label{step0}
\eea
where $K_n$ is a cut-off which increases with $n$ (chosen below) and $R_n$ 
is the remainder of 
the Taylor series for $\log(I-M_n)$ written in the integral form:
\bea
R_n(K) = \int_0^1 \tr\left(\frac{M_n^{K+1}}{(1-xM_n)^{K+1}}\right)
(1-x)^K dx. \nonumber
\eea
An upper bound on $|R_n(K)|$ follows from Lemma \ref{lm2} by replacing all eigenvalues
of $M_n$ with $\lambda_{max}(n)$:
\bea
|R_n(K)| & \leq &  n\lambda_{max}^{K+1}(n) \int_0^1 \frac{(1-x)^K }{(1-\lambda_{max}(n)x)^{K+1}}
dx \nonumber \\
& \leq &
 n\lambda_{max}^{K}(n) \log\left(\frac{1}{1-\lambda_{max}(n)}\right) \nonumber \\
& \leq & n\log\left(\frac{n}{\mu}\right)
\left(1-\frac{\mu}{n}\right)^{K}. \nonumber
\eea
So, if we choose
\bea\label{cutoff}
K_n=\left \lfloor n^\alpha \right \rfloor, ~\alpha>1,
\eea
it is easy to check that
\bea\label{rem}
\lim_{n\rightarrow \infty} R_n\left(K_n\right)=0.
\eea

The last step of the proof is the calculation of $\sum_{m=1}^{K_n}\frac{1}{m}\tr M_n^m$. 
The relevant results can be summarized as follows:
\begin{lemma}\label{lm3}
For any fixed integer $m>0$,
\bea \label{trace}
\lim_{n\rightarrow \infty} \frac{1}{\sqrt{2n}} \mathrm{Tr} M_n^m=\sqrt{\frac{1}{2\pi m}}.
\eea
Moreover, for any any positive integers  $m,n$
\bea\label{ub}
\mathrm{Tr} M_n^m\leq \sqrt{\frac{n}{\pi m}}\left(1+n^{-1}\right) +\frac{1}{4}
+\frac{1}{8} \sqrt{\frac{m}{\pi n}}\left(1+ 2n^{-1}\right).
\eea
\end{lemma}
Let us stress that formula (\ref{trace}) alone is not enough for the calculation of the $\lim_{n \rightarrow \infty} n^{-1/2} \log p_{2n,0}$ using (\ref{step0}) since the limits $n\rightarrow \infty$ and $m\rightarrow \infty$ do not necessarily commute. 
Instead, let us fix an arbitrary integer $K>0$. For a sufficiently large $n$ (so that $K_n>K$) relation (\ref{ub}) gives:
\bea
&&\frac{1}{\sqrt{2n}}\sum_{m=1}^{K} \frac{1}{m}\tr M_n^m\leq \frac{1}{\sqrt{2n}}\sum_{m=1}^{K_n} 
\frac{1}{m}\tr M_n^m \nonumber \\
&& \hspace{.2in} \leq \frac{(1+n^{-1})}{\sqrt{2\pi}}\sum_{m=1}^{K_n} m^{-3/2}+\frac{1}{4\sqrt{2n}} 
\sum_{k=1}^{K_n} \frac{1}{m}+\frac{(1+ 2n^{-1})}{8\sqrt{2\pi } \; n}\sum_{k=1}^{K_n}\frac{1}{\sqrt{m}}.
\label{din}
\eea
In writing the above double inequality we used the fact that $M_n$ is positive definite, which implies that
$\tr M_n^m>0$ for all values of $m,~n$. Let us choose $K_n$ in the form (\ref{cutoff}) with $\alpha<2$
and take $n\rightarrow \infty$   in (\ref{din}). As $K$ is $n$-independent, we can use formula (\ref{trace})
to compute the limit of the left hand side. On the right hand side, the last two sums vanish in the limit
(as $\log(n)/\sqrt{n}$ and $n^{\alpha/2-1}$ correspondingly). The first sum converges to
$$\frac{1}{\sqrt{2\pi}}\sum_{m=1}^{\infty} m^{-3/2}=\frac{1}{\sqrt{2\pi}}\zeta(3/2),$$
where $\zeta(x)=\sum_{m=1}^\infty m^{-x}$ is the Riemann zeta-function. 

We have found that for any positive integer $K$,
\bea
\frac{1}{\sqrt{2\pi}}\sum_{m=1}^{K}m^{-3/2} & \leq & \liminf_{n\rightarrow \infty}\frac{1}{\sqrt{2n}}\sum_{m=1}^{K_n} \frac{1}{m} \tr M_n^m
\nonumber \\
& \leq &  \limsup_{n\rightarrow \infty}\frac{1}{\sqrt{2n}}\sum_{m=1}^{K_n} \frac{1}{m} \tr M_n^m \nonumber \\
& \leq & \frac{1}{\sqrt{2\pi}}\zeta(3/2). \nonumber 
\eea
As $K$ is arbitrary, we conclude that
\bea\label{final}
\lim_{n\rightarrow \infty}\frac{1}{\sqrt{2n}}\sum_{m=1}^{K_n} \tr M_n^m= \frac{1}{\sqrt{2\pi}}\zeta(3/2).
\eea
So we proved that both (\ref{rem}) and  (\ref{final}) hold provided the cut-off is taken
in the form (\ref{cutoff}) for any fixed $\alpha \in (1,2)$.  

Finally, we can take the $n\rightarrow \infty$ limit in (\ref{step0}). Employing (\ref{rem}) and (\ref{final})
we find that
\bea 
\label{thm1k0}
\lim_{n\rightarrow \infty} \frac{1}{\sqrt{2n}}\log p_{2n,0}=-\frac{1}{\sqrt{2\pi}}\zeta(3/2). \nonumber
\eea
Theorem \ref{thm1} is proved for $k=0$.
\subsection{The calculation of $\lim_{n\rightarrow \infty} \frac{1}{\sqrt{2n}} \log p_{2n,2k}$ for $k>0$.}
It follows from Lemma \ref{lm1} that
\bea
p_{2n,2k}=\frac{1}{k!} \left(\frac{d}{dz}\right)^k\det(I+(z-1)M_n)\mid_{z=0}. \nonumber
\eea
Equivalently,
\bea 
\label{pnke}
p_{2n,2k}=\frac{p_{2n,0}}{k!} \left(\frac{d}{dz}\right)^k\det(I+zP_n)\mid_{z=0},
\eea
where $P_n=(I-M_n)^{-1}M_n$. Recall that
\[
\det(I+zP_n)=\sum_{k=0}^nz^k e_{k}(\nu),
\]
where $\nu=(\nu_1, \nu_2, \ldots, \nu_n)$ are the eigenvalues of $P_n$ and
$e_{k}$ is the degree-$k$ elementary symmetric polynomial in $n$ variables \cite{mcd},
\bea
e_{k}(\nu)=
\sum_{1\leq i_1<i_2<\ldots <i_k\leq n}\nu_{i_{1}}\nu_{i_{2}}\ldots \nu_{i_{k}}. \nonumber
\eea 
Therefore,
\bea\label{symfun}
p_{2n,2k}=p_{2n,0} \, e_{k}(\nu),  \quad \mbox{for $k=0,1,\ldots, n.$}
\eea
Let us enumerate the eigenvalues of $M_n$ and $P_n$ as follows:
$$\lambda_1\geq\lambda_2\geq\ldots \geq\lambda_n>0,$$
$$\nu_1\geq\nu_2\geq\ldots \geq\nu_n>0.$$
By the definition of $P_n$, $\nu_i=\frac{\lambda_i}{1-\lambda_i}$.
Note that $\nu_i$ is
a monotonically increasing function of $\lambda_i$.
Combining this remark with the spectral bound of Lemma \ref{lm2}, we get the following bound on the elementary
symmetric polynomials:
\bea\label{intstep}
e_{k}(\nu)\leq \nu_1^k e_{k}(1,1,\ldots,1)\leq
\left(\frac{\lambda_1}{1-\lambda_{1}}\right)^kn^k\leq 
\left(\frac{n}{\mu}\right)^k n^k.
\eea
Substituting (\ref{intstep}) into (\ref{symfun}) we obtain the following upper bound on $\log p_{2n,2k}$:
\bea\label{part1}
\log p_{2n,2k} \leq \log p_{2n,0}+k\log\left(\frac{n^2}{\mu}\right).
\eea
Next we derive a lower bound on $\log p_{2n,2k}$. 
By positive definiteness, $\nu_i\geq \lambda_i$ and therefore $e_{k}(\nu)\geq  e_{k}(\lambda)$.
Let us fix a positive integer $k$.
Due to (\ref{trace}), for any $\epsilon>0$ there is a positive integer
$N_\epsilon$ such that
for any $n>N_\epsilon$

\bea \label{rhobnd}
 \sqrt{\frac{n}{\pi }}(1-\epsilon)\leq
\tr M_n \leq \sqrt{\frac{n}{\pi }}(1+\epsilon).
\eea
On the other hand,
\begin{eqnarray}\label{labnd}
\tr M_n=(\lambda_1+\ldots+\lambda_{k-1})
+(\lambda_k+\ldots+\lambda_n)\leq (k-1)+(n-k+1)\lambda_k, 
\end{eqnarray}
where the inequality is due to (\ref{specrad}) and the chosen ordering of
$\lambda$'s.

Combining (\ref{rhobnd}) and (\ref{labnd})
we obtain the following bound on the $k$-th largest eigenvalue
of $M_n$:

\bea\label{labnd2}
\lambda_{k}\geq \frac{\sqrt{\frac{n}{\pi}}(1-\epsilon)-k+1}{n-k+1},
\eea
which holds for $n>N_\epsilon$.
Inequality (\ref{labnd2}) leads to the desired bound for $e_k(\nu)$:
\bea
e_k(\nu)\geq e_k(\lambda)\geq \lambda_1\lambda_2\ldots \lambda_k
\geq \lambda_k^k\geq \left(\frac{\sqrt{\frac{n}{\pi}}(1-\epsilon)-k+1}{n-k+1}\right)^k. \nonumber
\eea
Substituting this result into (\ref{symfun}) we find that 
\bea\label{part2}
\log p_{2n,2k} \geq \log p_{2n,0}+k\log\left(\frac{\sqrt{\frac{n}{\pi}}(1-\epsilon)-k+1}{n-k+1}\right).
\eea

Combining (\ref{part1}) and (\ref{part2}) we find that
\bea\label{limklim0}
\lim_{n\rightarrow \infty} \frac{1}{\sqrt{2n}} \log p_{2n,2k}=
\lim_{n\rightarrow \infty} \frac{1}{\sqrt{2n}} \log p_{2n,0}.
\eea
Relations (\ref{limklim0}) and (\ref{thm1k0}) imply that formula (\ref{result})
of Theorem \ref{thm1}
is proved for any fixed integer $k>0$.
 
Moreover, it is evident from (\ref{part1}) and (\ref{part2}) that the equality
(\ref{limklim0}) generalizes to
\bea\label{limknlim0}
\lim_{n\rightarrow \infty} \frac{1}{\sqrt{2n}} \log p_{2n,2k_n}=
\lim_{n\rightarrow \infty} \frac{1}{\sqrt{2n}} \log p_{2n,0}.
\eea
where $(k_n)_{n\geq 1}$ is a sequence of natural numbers such that
$$\lim_{n\rightarrow \infty} \left(n^{-1/2}\log(n)\right) k_{n}=0.$$ This proves the last
claim of Theorem \ref{thm1}.\\
\\
{\bf Remark}. Our proof of the $k>0$ part of the Theorem 
is a simple consequence of positive-definiteness
of $M_n$, the spectral bound and the fact that $\tr(M_n)\stackrel{n\rightarrow \infty}{\simb} \sqrt{n/\pi}$. It is interesting that the proof does not rely on any detailed knowledge of the spectrum of $M_n$.
\section*{Acknowledgements}
C.T. acknowledges useful discussions with K. Nestmann.
This research was supported by the Israel Science Foundation through the grant 
No.\ 647/12 (E.K.); an EPSRC grant No.\ EP/K011758/1 (M.P. and R.T.); and  a Leverhulme Trust Research Fellowship (O.Z.).

\begin{appendix}
\section{Proofs for the lemmas.}\label{app1}
\subsection{Lemma 1.}
To prove the lemma, we start with the exact formula due to Kanzieper and Akemann \cite{KA-2005} which expresses the probabilities $p_{2n,2k}$
in terms of elementary symmetric functions:
\begin{eqnarray}\label{KA-solution-1}
    p_{2n,2k} = p_{2n,2n} \, e_{n-k} (t_1,\dots,t_{n-k}),
\end{eqnarray}
where $t_j$'s are given by
\begin{eqnarray} \label{KA-solution-2}
    t_j = \frac{1}{2} {\rm \tr}{({\bm A}^{-1}{\bm B})}^j.
\end{eqnarray}
Here ${\bm A}$ and ${\bm B}$ are $2n\times 2n$ antisymmetric matrices whose entries
\begin{eqnarray} \label{A-m}
    {\bm A}_{jk} &=& \left<q_{j-1},q_{k-1}\right>_{{\rm R}},\\
    \label{B-m}
    {\bm B}_{jk} &=& \left<q_{j-1},q_{k-1}\right>_{{\rm C}},
\end{eqnarray}
are defined in terms of skew products
\begin{eqnarray}\label{sops-1}
    \left<f,g\right>_{{\rm R}} =
    \frac{1}{2} \int_{{\mathbb R}^2} dx \,dy\, e^{-(x^2+y^2)/2} {\rm sgn}(y-x) \, f(x)\, g(y)
\end{eqnarray}
and
\begin{eqnarray}
   \left<f,g\right>_{{\rm C}} = i \int_{{\rm Im\,}z >0} d^2z \,  e^{-(z^2+\bar{z}^2)/2} \, {\rm erfc} \left(
    \frac{z-\bar{z}}{i\sqrt{2}}
    \right)
     \, \left[ f(z)\, g(\bar{z}) - g(z)\, f(\bar{z}) \right].
\end{eqnarray}
Let us stress that (\ref{KA-solution-1}) is valid for an 
{\it arbitrary} choice of monic polynomials $q_j(x)$ of degree $j$, provided  matrix ${\bm A}$ is invertible.
\newline\newline\noindent
Substituting Eqs.~(\ref{KA-solution-1}) and (\ref{KA-solution-2}) into the generating function
\begin{eqnarray}
    g_{2n}(z) = \sum_{k=0}^n z^k p_{2n,2k}
\end{eqnarray}
and making use of the summation formula \cite{mcd}
\begin{eqnarray}
    \sum_{\ell=0}^\infty z^\ell e_\ell (t_1,\dots,t_\ell) = \exp\left(
        \sum_{j=1}^\infty (-1)^{j-1} t_j \,\frac{z^j}{j}
    \right),
\end{eqnarray}
we obtain the Pfaffian representation \cite{KA-2005,eugene1,BK-2007}:
\begin{eqnarray} \label{e1}
    g_{2n}(z) = p_{2n,2n} \,{\rm Pf} (-{\bm A}^{-1})\, {\rm Pf} (z {\bm A} + {\bm B}),
\end{eqnarray}
see remark $1.3$ of \cite{BK-2007} justifying the transition from  square
roots of determinants to Pfaffians.
Since $g_{2n}(1)=1$, $p_{2n,2n}=\left({\rm Pf} (-{\bm A}^{-1})\, {\rm Pf} (z {\bm A} + {\bm B})\right)^{-1}$ and (\ref{e1}) simplifies to
\begin{eqnarray} \label{e4}
g_{2n}(z) = \frac{{\rm Pf}(z{\bm A}+{\bm B})}{{\rm Pf}({\bm A}+{\bm B})}.
\end{eqnarray}
Next we will use the fact that expression (\ref{e4}) for the generating function 
does not depend on a particular choice of monic polynomials $q_j(x)$ in (\ref{A-m}) and (\ref{B-m}) to simplify it even further. 
Namely, we will choose $q_j(x)$'s in such a way that the matrix ${\bm A}+{\bm B}$ is block diagonal. Clearly, such polynomials should be skew-orthogonal with respect to the skew product
\begin{eqnarray}
    \left<f,g\right> = \left<f,g\right>_{{\rm R}}+ \left<f,g\right>_{{\rm C}},
\end{eqnarray}
that is
\begin{eqnarray}
    \left<
    q_{2j}, q_{2k+1}
    \right> = - \left<
    q_{2k+1}, q_{2j}
    \right> = {r}_j \delta_{j,k},\qquad \left<
    q_{2j}, q_{2k}
    \right> = \left<
    q_{2j+1}, q_{2k+1}
    \right> = 0.
\end{eqnarray}
These were first calculated in the paper \cite{forrester_nagao}:
\begin{eqnarray}\label{rcoeff}
    q_{2j}(x) = x^{2j},\qquad q_{2j+1}(x) = x^{2j+1} - 2j x^{2j-1}, \qquad {r}_j = \sqrt{2\pi} \Gamma(2j+1).
\end{eqnarray}
Given the choice of $q_j$'s described above,\newline
\newline\noindent
(a) the matrix ${\bm A}+{\bm B}$ acquires a block-diagonal form
\begin{eqnarray}
 {\bm A}+{\bm B} = {\bm r} \otimes {\bm J},\qquad {\bm r} = {\rm diag} \left({r}_0,\dots,{r}_{n-1}\right), \qquad
 {\bm J} = \left(
                           \begin{array}{cc}
                             0 & 1 \\
                             -1 & 0 \\
                          \end{array}
                         \right),
\end{eqnarray}
which leads to
\begin{eqnarray} \label{e5}
g_{2n}(z) = \frac{{\rm Pf}({\bm r} \otimes {\bm J} + (z-1){\bm A})}{{\rm Pf}({\bm r} \otimes {\bm J})}.
\end{eqnarray}
\noindent
(b) the matrix ${\bm A}$ is given by
\begin{eqnarray}\label{acoeff}
    {\bm A}_{2j,2k} = {\bm A}_{2j+1,2k+1} = 0, \qquad {\bm A}_{2j-1,2k} = \Gamma\left(
    j+k-\frac{3}{2}
    \right).
\end{eqnarray}
Notice that matrix elements of both 
${\bm r} \otimes {\bm J}$ and ${\bm A}$ labeled by a pair of  indexes of the same parity vanish. Therefore, the  $2n\times 2n$ Pfaffians in the numerator
and the denominator of (\ref{e5}) are reduced to $n\times n$ determinants: 
\begin{eqnarray} \label{e44}
g_{2n}(z) = \frac{{\rm det}[r_{j-1}\delta_{jk} + (z-1){\bm A}_{2j-1,2k}]_{1\leq j,k\leq n}}{{\rm det}[r_{j-1}\delta_{jk}]_{1\leq j,k\leq n}}.
\end{eqnarray}
Finally, we apply the formula $\det(U)/\det(V^2)=\det(V^{-1}UV^{-1})$ to perform division in (\ref{e44}). With the help of the explicit formulae (\ref{rcoeff})
and (\ref{acoeff}) we get
\begin{eqnarray} \label{e45}
g_{2n}(z) = {\rm det}\left[\delta_{jk} + \frac{(z-1)}{\sqrt{2\pi}}
\frac{\Gamma(j+k-\frac{3}{2})}{\sqrt{\Gamma(2j-1)\Gamma(2k-1)}}\right]_{1\leq j,k\leq n}.
\end{eqnarray}
Lemma \ref{lm1} is proved.

\subsection{Lemma 2.}
The proofs of Lemmas 2,~3 are based on the following integral representation for the matrix elements
 (\ref{mat})
of matrix $M_n$:
\bea\label{matint}
M_n(j,k)=\frac{1}{\sqrt{2\pi}}\int_{0}^\infty \frac{dx}{x^{5/2}}e^{-x}
\frac{x^j}{\sqrt{\Gamma(2j-1)}}
\frac{x^k}{\sqrt{\Gamma(2k-1)}},~1\leq j,k\leq n,
\eea
which can be obtained by representing $\Gamma(j+k-3/2)$ in (\ref{mat}) as an integral.

Take any $v=(v_1, v_2, \ldots, v_n) \in \R^n\setminus \{0\}$.
It follows from (\ref{matint}) that
\bea
\langle v,M_n v\rangle= \frac{1}{\sqrt{2\pi}}
\int_{0}^\infty \frac{dx}{x^{5/2}}e^{-x}
\left(\sum_{j=1}^n\frac{v_jx^j}{\sqrt{\Gamma(2j-1)}}\right)^2>0.
\eea
So, $M_n$ is positive definite by definition. 

Next, let us prove bound (\ref{specrad}) on the spectral radius of $M_n$. 
Let $\lambda_1, \lambda_2, \ldots, \lambda_n>0$ be the eigenvalues of $M_n$. Then
\bea\label{lm3proof}
\lambda_{max}(n)=\left(\lambda_{max}^n(n)\right)^{\frac{1}{n}}
\leq \left(\sum_{k=1}^n \lambda_k^n\right)^{\frac{1}{n}}
=\left(\tr M_n^n \right)^{\frac{1}{n}}.
\eea
It follows from the upper bound (\ref{ub}) of Lemma \ref{lm3} that
for any $\epsilon>0$, there is $N_\epsilon$ such that for any $n>N_{\epsilon}$,
\[
\tr M_n^n \leq \sqrt{\frac{1}{\pi}}+ \frac14 +\frac{1}{8}\sqrt{\frac{1}{\pi}}
+ \epsilon=1-c_1+ \epsilon,
\]
where $c_1\in (0,1)$. Therefore, we can choose $\epsilon$ small enough so that 
\[
\tr M_n^n \leq 1-\mu,
\]
where $\mu \in (0,1)$. Using this estimate in (\ref{lm3proof}) for $n>N_\epsilon$
we get
\bea
\lambda_{max}(n)\leq (1-\mu)^{\frac{1}{n}}\leq 1-\frac{\mu}{n}.
\eea
Lemma \ref{lm2} is proved  for $N=N_\epsilon$.\\
\\
{\bf Remark.} 
The spectral properties of $M_n$ seem quite interesting. For instance, in the large-$n$ limit there is
a concentration of eigenvalues near $1$ 
such that the restriction of $M_n$ to the space spanned by the corresponding eigenvectors
is close to an identity operator perturbed by an elliptic linear differential
operator. Formal analysis of this perturbation suggests the asymptotic
$\lambda_{max}(n) = 1 - \mu_0 n^{-1} + \underline{o}(n^{-1})$ for suitable $\mu_0 >0$. 
\subsection{Lemma 3.} 
The integral representation (\ref{matint}) for the matrix elements of $M_n$ leads
to the following integral representation for the trace of a power of $M_n$:

\bea\label{ir1}
\tr M_n^m&=&\int_{0}^\infty \frac{dx_1}{\sqrt{2\pi x_1}} \int_{0}^\infty \frac{dx_2}{\sqrt{2\pi x_2}}
\ldots \int_{0}^\infty \frac{dx_m}{\sqrt{2\pi x_m}} e^{-x_1-x_2-\ldots - x_m}\nonumber\\
&&\cosh_{n-1}(\sqrt{x_{m} x_1})\cosh_{n-1}(\sqrt{x_1x_2})\ldots \cosh_{n-1}(\sqrt{x_{m-1} x_m}),
\eea
where
$
\cosh_{n}(x)=\sum_{k=0}^n\frac{x^{2k}}{(2k)!}
$
is the degree-$2n$ Taylor polynomial generated by the hyperbolic cosine. 
Performing the change of variables $x_k=y_k^2$ in (\ref{ir1})
we can re-write the integral representation for $\tr M_n^m$ as follows:
\bea\label{ir2}
\tr M_n^m&=&\left(\frac{2}{\pi}\right)^{m/2}\int_{\R_{+}^m} dy
e^{-\sum_{k=1}^m y_k^2}\nonumber\\
&&\cosh_{n-1}(y_{m} y_1) \cosh_{n-1}(y_1y_2)\ldots \cosh_{n-1}(y_{m-1} y_m),
\eea
Here $\R_{+}^m=\{(y_1, y_2,\ldots, y_m)\in \R^m\mid y_k\geq 0,~ k=1,2,\ldots, m\}$ is the first
`quadrant' of $\R^m$ and $dy$ is a shorthand notation for Lebesgue measure on $\R^m$.
As the integrand of (\ref{ir2}) is symmetric with respect to reflection $y_i\rightarrow -y_i$
for any $i=1,2,\ldots, m$, we can re-write $\tr M_n^m$ as an integral over $\R^m$:
\bea\label{ir3}
\tr M_n^m&=&\left(\frac{1}{2\pi}\right)^{m/2}\int_{\R^m} dy
e^{-\sum_{k=1}^m y_k^2}\nonumber\\
&&\cosh_{n-1}(y_{m} y_1) \cosh_{n-1}(y_1y_2)\ldots \cosh_{n-1}(y_{m-1} y_m).
\eea

To prove Lemma 3 we will establish an upper and a lower bound on $\tr M_n^m$
and then compute the large-$n$ limit of each of these bounds.
\subsubsection{An upper bound for $\tr M_n^m$.}
A good starting point for the calculation is formula (\ref{ir3}). 
For any $x\in \R$, $\cosh_{n-1}(x)\leq \cosh(x)$. Also,
\bea\label{ci}
\cosh_{n-1} (x)=\oint \frac{dz}{2\pi i z} \frac{1-z^{-2n}}{1-z^{-2}}e^{zx}, 
\eea
where the integral is anti-clockwise around a circle of radius smaller than $1$ centred at the origin in the complex plane.
Replacing all but one $\cosh_{n-1}$ with $\cosh$ we get:
\bea\label{upper1}
\tr M_n^m&\leq& \left(\frac{1}{2\pi}\right)^{m/2}\int_{\R^m} dy
e^{-\sum_{k=1}^m y_k^2}
\cosh_{n-1}(y_{m} y_1) \cosh(y_1y_2)\ldots \cosh(y_{m-1} y_m)\nonumber\\
& = & \left(\frac{1}{2\pi}\right)^{m/2}\mathbb{E}_{\alpha_1 \alpha_2\ldots \alpha_{m-1}}\int_{\R^m} dy
e^{-\sum_{k=1}^m y_k^2}
\cosh_{n-1}(y_{m} y_1) e^{\sum_{l=1}^{m-1} \alpha_l y_{l}y_{l+1}},
\eea
where $\alpha_1, \alpha_2, \ldots, \alpha_{m-1}$ are independent identically distributed random variables which take values $\pm 1$ with probability $1/2$.
Representing the remaining $\cosh_{n-1}$ with the help of (\ref{ci}) and then computing resulting 
Gaussian integral over $\R^m$ we find
\bea\label{upper2}
\tr M_n^m\leq \left(\frac{1}{2}\right)^{m/2} \mathbb{E}_{\alpha_1 \alpha_2\ldots \alpha_{m-1}}
\oint \frac{dz}{2\pi i z} \frac{1-z^{-2n}}{1-z^{-2}}[D_m^{(\alpha)}(z)]^{-1/2}, 
\eea  
where
\bea\label{mdet}
D_m^{(\alpha)}(z)=\det\left(
\begin{array}{ccccccc}
1 & -\frac{\alpha_1}{2} & 0 & 0 &\ldots &0&-\frac{z}{2} \\
-\frac{\alpha_1}{2} & 1 & -\frac{\alpha_2}{2}& 0&0&\ldots &0 \\
0 &-\frac{\alpha_2}{2}&1&-\frac{\alpha_3}{2}&0 & \ldots & 0 \\
\cdot& \cdot& \cdot& \cdot& \cdot& \cdot& \cdot\\
\cdot& \cdot& \cdot& \cdot& \cdot& \cdot& \cdot\\
\cdot& \cdot& \cdot& \cdot& \cdot& \cdot& \cdot\\
0& \ldots& 0& -\frac{\alpha_{m-3}}{2}& 1& -\frac{\alpha_{m-2}}{2}& 0\\
0& \ldots& 0& 0& -\frac{\alpha_{m-2}}{2}& 1& -\frac{\alpha_{m-1}}{2}\\
-\frac{z}{2}& 0& \ldots& 0& 0& -\frac{\alpha_{m-1}}{2}& 1
\end{array} \right).
\eea
The determinant can be calculated recursively in $m$, yielding $D_1^{(\alpha)}(z) = 1-z$ and
\bea\label{det}
D_m^{(\alpha)}(z)=-(m-1)\frac{1}{2^m} \left(z-A_m\right)\left(z+A_m\frac{m+1}{m-1}\right) \quad 
\mbox{for $m \geq 2$,}
\eea
where $A_m=\prod_{k=1}^{m-1}\alpha_k$. 
Note that (\ref{det}) implies that all principal minors of the matrix under the sign
of the determinant in (\ref{mdet}) are positive for $z=0$. Therefore the matrix itself  
is positive definite for $z=0$. By continuity, the real part of
this matrix remains positive definite for $z \neq 0$ provided $|z|$ is small enough. 
Therefore, the real part of the quadratic form which determines the Gaussian integral in (\ref{upper1})
is positive definite, which justifies the interchange of integrals leading to (\ref{upper2}) provided the 
contour is taken to be a circle around the origin of a sufficiently small radius.

Substituting (\ref{det}) into (\ref{upper2})
and changing the integration variable $z\rightarrow A_m z$
 we find that the integrand no longer depends on $\alpha$'s. Averaging over $\alpha$'s
becomes trivial and we get the following
integral upper bound
\bea\label{upper3}
\tr M_n^m\leq \oint \frac{dz}{2\pi z} \frac{z^{-2n}-1}{z^{-2}-1} \frac{1}{\sqrt{1-z}}\frac{1}{\sqrt{(m-1)z+m+1}}.
\eea
The rest of the calculation is slightly different depending on whether $m=1$ or $m>1$.
Here present the calculation for $m>1$ only, the (simpler) case of $m=1$ can be treated
along similar lines.
We calculate the integral in the right hand side of (\ref{upper3}) as follows. First we
replace $z^{-2n}-1$ with $z^{-2n}$ in the integrand on the r.h.s. of (\ref{upper3}), since 
this does not change the value of the integral as the omitted term is analytic inside of the contour
of integration. Next we deform the contour away from the singularity at zero and out to infinity, leading to integrals around the other
singularities of the (modified) integrand: a simple pole at $z=-1$, a branch cut singularity along the real line from $1$ to $+\infty$, and
a branch cut singularity along the real line from $-\frac{m+1}{m-1}$ to $-\infty$. The contribution from the 
the integral over the large circle at infinity is zero. The contribution from the pole at $z=-1$ is easily evaluated as
$1/4$.  Evaluating the integral around the branch from $1$ to $+\infty$ it is convenient first integrate by parts, so that 
the singularity at $z=1$ is integrable. The integrals along the two branch cuts lead to two real integrals 
whose asymptotics are controlled by the integrand $(1+y)^{-2n}$. Changing variable $y \to y/2n$,
and making some simple estimates on terms that do not affect the leading asymptotics, we are led to 
\bea
&\tr M_n^m& \leq \frac{1}{4}+\sqrt{\frac{n}{\pi m}}
\int_{0}^\infty \frac{dy}{\sqrt{\pi y}}
\left(1+\frac{y}{2n}\right)^{-2n} \label{123456} \\
&&+\frac{1}{\sqrt{2\pi n}}\frac{m+1}{2\sqrt{m-1}}\left(\frac{m+1}{2m}\right)^{3/2}
\left(\frac{m-1}{m+1}\right)^{2n+1}\int_{0}^\infty \frac{dy}{\sqrt{\pi y}}
\left(1+\frac{y}{2n}\right)^{-2n+1}.\nonumber
\eea
Both integrals in the above expression can be estimated using the following bound:
\bea
I_M = \int_0^\infty \frac{dy}{\sqrt{\pi y}}
\left(1+\frac{y}{M}\right)^{-M}\leq 1+ \frac{2}{M},
\eea
which follows by evaluating the integral, using the substitution $t=(1+\frac{y}{M})^{-1}$, in terms of the beta function as 
\[
I_M = \sqrt{\frac{M}{\pi}} B\left(M-\frac12,\frac12\right) = \sqrt{M} \Gamma(M-\frac32)/\Gamma(M)
\] 
and using bounds on the Gamma function. 
Using this in (\ref{123456}), the final result is
\bea
\tr M_n^m & \leq & \frac14 + \sqrt{\frac{n}{\pi m}} \left(1+ \frac{1}{n}\right) 
+\frac{1}{8}\sqrt{\frac{m}{\pi n}} \left(1+\frac{2}{n}\right) \left(\frac{m-1}{m+1}\right)^{2n-\frac32} \nonumber \\
& \leq & \frac14 + \sqrt{\frac{n}{\pi m}} \left(1+ \frac{1}{n}\right) 
+\frac{1}{8}\sqrt{\frac{m}{\pi n}} \left(1+\frac{2}{n}\right) \label{upper4}
\eea
which coincides with the claim (\ref{ub}) of Lemma \ref{lm3}.

Dividing both sides of (\ref{upper4}) by $\sqrt{2n}$ and taking the large $n$ limit,
we find that
\bea\label{half_limit}
\limsup_{n\rightarrow \infty} \frac{1}{\sqrt{2n}}\tr M_n^m\leq \sqrt{\frac{1}{2\pi m}}
\eea
\subsubsection{The limit $\lim_{n\rightarrow \infty}\frac{1}{\sqrt{2n}}\tr M_n^m$.}
The strategy is to derive an integral lower bound for $\tr M_n^m$ and calculate
the large $n$-limit of the bound. Our starting point is the relation (\ref{ir2}) and
the following estimate for the polynomial $\cosh_{n-1}$:
\begin{lemma}\label{prop1}
There exist  two sequences $(h_n)_{n\geq 1}, (S_n)_{n\geq 1} \subset \R$
such that
\bea
&&\lim_{n\rightarrow \infty} h_n=\frac{1}{2},~\lim_{n\rightarrow \infty} S_n=2, \nonumber \\
&&e^{-ny}\cosh_{n-1}(ny)\geq h_n \mathbbm{1} (y< S_n), \quad \mbox{for $y\geq 0, ~ n\geq 1$.}
\label{lb0}
\eea
Here $\mathbbm{1}(y < S_n)$ is the indicator function of the set $[0, S_n)$.
\end{lemma}
In fact, as $n\rightarrow \infty$, $e^{-ny}\cosh_{n-1}(ny)$ converges almost everywhere to
$\frac{1}{2}\mathbbm{1}(y<2)$ for $y \geq 0$, but here
we only need the lower bound. The proof of Lemma \ref{prop1} is given 
in Section \ref{proofpr1}.

Using the bound (\ref{lb0}) in (\ref{ir2}) we find that
\bea\label{lb1}
\tr M_n^m\geq h_n^m \left(\frac{2}{\pi}\right)^{m/2} n^{m/2}
\int_{\R_{+}^m} dy \prod_{l=1}^m \mathbbm{1} (y_ly_{l+1}<S_n)
e^{-\frac{n}{2}\sum_{k=1}^m (y_{k+1}-y_k)^2},
\eea
where $y_{m+1}:= y_1$. 
It is straightforward to verify that the domain of integration
for the integral in (\ref{lb1}) contains the hypercube $(0,\sqrt{S_n})^m$,
$$(0,\sqrt{S_n})^m\subset
\{y\in R_{+}^m | y_{k}y_{k+1} <S_n, k=1,2,\ldots,m\}.$$
Therefore,
\bea\label{rem1}
\prod_{l=1}^n \mathbbm{1} \left(y_l<\sqrt{S_n}\right)\leq \prod_{l=1}^m \mathbbm{1} (y_ly_{l+1}<S_n), ~y\in \R^m_{+}
\eea
Substituting (\ref{rem1}) in (\ref{lb1})  and changing the integration variables
according to
\bea\label{chofv}
R&=&y_1+y_2+\ldots+y_m, \nonumber\\
z_k&=&y_{k+1}-y_k,~k=1,2,\ldots, m-1, \nonumber
\eea
we get the following lower bound:
\bea\label{lb2}
&& \hspace{-.3in} \tr M_n^m \geq  \frac{ h_n^m}{m} \left(\frac{2}{\pi}\right)^{m/2} n^{m/2}
\int_{0}^{m\sqrt{S_n}}dR  \nonumber \\
&& \hspace{.5in} \int_{P_{m-1}(R)}dz_1\ldots dz_{m-1} 
e^{-\frac{n}{2}[\sum_{k=1}^{m-1} z_k^2+(\sum_{k=1}^{m-1} z_k)^2]},
\eea
where $P_{m-1}(R)$ is the intersection of the hypercube $(0,\sqrt{S_n})^m$ and
the hyperplane $$\{y\in R_{+}^m|y_1+y_2+\ldots+y_{m}=R\}.$$
In the derivation of (\ref{lb2}) we used the fact that the Jacobian of the transformation $y\rightarrow (R,z)$ is equal to $1/m$. 

The large-$n$ limit of the right
hand side of (\ref{lb2}) can be evaluated by arguing as in the Laplace method:
\begin{eqnarray*}
&& \liminf_{n\rightarrow \infty}\frac{1}{\sqrt{2n}} \tr M_n^m\\
& \geq & \lim_{n\rightarrow \infty} \frac{ h_n^m}{\sqrt{2n}m} \left(\frac{2}{\pi}\right)^{m/2} n^{m/2}
\int_{0}^{m\sqrt{S_n}}dR \int_{\R^{m-1}}dz_1\ldots dz_{m-1} 
e^{-\frac{n}{2}[\sum_{k=1}^{m-1} z_k^2+(\sum_{k=1}^{m-1} z_k)^2]}\\
& = & \lim_{n\rightarrow \infty} \sqrt{\frac{S_n}{2}} h_n^m  \left(\frac{2}{\pi}\right)^{m/2} n^{(m-1)/2}
\int_{-\infty}^\infty\frac{d\lambda}{2\pi}
\int_{\R^{m}} dz_1 \ldots dz_{m} e^{i\lambda\sum_{k=1}^m z_k}
e^{-\frac{n}{2}\sum_{k=1}^m z_k^2}\\
& = & \lim_{n\rightarrow \infty} \sqrt{\frac{S_n}{2}} h_n^m  \left(\frac{2}{\pi}\right)^{m/2} n^{(m-1)/2}
\int_{-\infty}^\infty\frac{d\lambda}{2\pi}
\left(\int_{-\infty}^{\infty}dz e^{i\lambda z-\frac{n}{2}z^2}\right)^m\\
& =& \lim_{n\rightarrow \infty} \sqrt{\frac{S_n}{2}} h_n^m  \left(\frac{2}{\pi}\right)^{m/2} n^{(m-1)/2} \left(\frac{2\pi}{n}\right)^{\frac{m}{2}}
\int_{-\infty}^\infty\frac{d\lambda}{2\pi}
e^{-\frac{m}{2n}\lambda^2}\\
& = & \lim_{n\rightarrow \infty} \sqrt{\frac{S_n}{2}} h_n^m  \left(\frac{2}{\pi}\right)^{m/2} n^{(m-1)/2} \left(\frac{2\pi}{n}\right)^{\frac{m}{2}}
\sqrt{\frac{ n}{2\pi m}}=\sqrt{\frac{1}{2\pi m}}.
\end{eqnarray*}
The crucial, albeit very standard, first step in the above derivation consists of verifying that
extending the integration
space for the $z$-integral from $P_{m-1}(R)$, when $R\in (0,2)$, to $\R^{m-1}$  
doesn't change the large $n$-limit.

We conclude that 
\bea
\liminf_{n\rightarrow \infty}\frac{1}{\sqrt{2n}}\tr M_n^m \geq \sqrt{\frac{1}{2\pi m}}, \nonumber
\eea
and in combination with (\ref{half_limit}) this gives
\bea
\lim_{n\rightarrow \infty}\frac{1}{\sqrt{2n}}\tr M_n^m = \sqrt{\frac{1}{2\pi m}}. \nonumber
\eea
Statement (\ref{trace}) of Lemma \ref{lm3} is proved.
\subsection{Lemma 4.}\label{proofpr1}
Let $\left\{\alpha_n\right\}_{n=1}^{\infty}$ be an 
arbitrary sequence of positive real numbers which diverges as $n\to\infty$ slower than $n^{1/2}$, that is
$\lim\limits_{n\rightarrow \infty } \alpha_n=\infty$, but
	$\lim\limits_{n\to\infty} \alpha_n n^{-1/2}=0$. We will show that
	there exists $N_0>0$ such that for any $n>N_0$ and $x\geq 0$
	\begin{equation}\label{pp1}
		e^{-nx}\cosh_n\left(nx\right)\geq 
		\left(\frac{1}{2}-\frac{1}{\sqrt{4\pi}}
		\alpha_n^{-1}e^{-\alpha_n^2/4}\right)
		\mathbbm{1}(x\leq 2-\alpha_nn^{-1/2}).
	\end{equation}
The statement of Lemma \ref{prop1}, where $\cosh_n(nx)$ is replaced by $\cosh_{n-1}(nx)$, 
is easily deduced from equation  (\ref{pp1}).

Our proof builds on the ideas of \cite{BM:06} dedicated to the  
study of sections of exponential series (Taylor polynomials generated by $\exp$). 
Let $e_n$ be a section of 
	exponential series defined by 	
		\begin{equation*}
			e_n\left(x\right)=\sumd_{j=0}^{n} \dfrac{x^{j}}{j!}.
		\end{equation*}
Consider also
		\begin{equation*}
			\expp{n}{x}=e^{-nx}e_n\left(nx\right), \qquad
			\expm{n}{x}=e^{-nx}e_n\left(-nx\right).
		\end{equation*}	
	Then the function we are interested in can be written 
	as
		\begin{equation*}
			f_n\left(x\right):=e^{-nx}\cosh_n\left(nx\right)
			= \dfrac{1}{2}\left(\expp{2n}{\dfrac{x}{2}}+
			\expm{2n}{\dfrac{x}{2}}\right).
		\end{equation*}		
	First we show that $\expm{2n}{x}>0$ for $x\geq 0$. One can 
	check that
		\begin{equation}\label{neg}
			\der{x} \left(e^{2nx}e_{2n}\left(-2nx\right)\right)
			=\dfrac{1}{\left(2n-1\right)!}
			\left(2nx\right)^{2n}e^{2nx}\geq 0,
		\end{equation}	
	and $\left.e^{2nx}e_{2n}\left(-2nx\right)\right|_{x=0}=1$.
	So $\expm{2n}{x}\geq e^{-4nx}>0$ for $x\geq 0$. The next step
	is to show that $f_n\left(x\right)$ is a decreasing function.
	However
		\begin{equation*}
			f_n'\left(x\right)=-n\expm{2n}{\dfrac{x}{2}},
		\end{equation*}
	which is negative by (\ref{neg}). 

The fact that  $f_n\left(x\right)$ is decreasing and the positivity of $\expm{2n}{x}$ imply that for any non-negative $x$
		\begin{eqnarray*}
			f_n\left(x\right)&\geq& f_n\left(x\right)
			\mathbbm{1}(x\leq 2-\alpha_nn^{-1/2})
			\\ &\geq& f_n\left(2-\alpha_nn^{-1/2}\right)
\mathbbm{1}(x\leq 2-\alpha_nn^{-1/2})
			\\ &\geq& \dfrac{1}{2} 
			\expp{2n}{1-\dfrac{\alpha_n}{2}n^{-1/2}}
\mathbbm{1}(x\leq 2-\alpha_nn^{-1/2}).
		\end{eqnarray*}	
	Therefore, it remains to prove that
		\begin{equation}\label{desbnd}
			\expp{2n}{1-\dfrac{\alpha_n}{2}n^{-1/2}}
			\geq 1-\sqrt{\dfrac{1}{\pi}}
			\alpha_n^{-1}e^{-\alpha_n^2/4},
		\end{equation}	
	for all $n>N_0$, where $N_0$ is chosen to satisfy 
	$\alpha_nn^{-1/2}<2$ for all $n>N_0$.
		
	We start with a differential equation
	satisfied by $e_{n}^{(+)}$. As it is easy to check,
		\begin{equation}\label{e:d_expp}
			\der{x}\expp{n}{x} = -\dfrac{1}{\left(n-1\right)!}
			\left(nx\right)^ne^{-nx}.
		\end{equation}
	So $\expp{n}{x}$ is a decreasing function on $\R_+$.
	
Equation (\ref{e:d_expp}) has to 
be solved with a boundary condition $\lim_{x\rightarrow \infty}\expp{n}{x} = 0$,
which follows from the definition of $e_{n}^{(+)}$. The solution is
		\begin{equation}\label{eppsol}
			\expp{n}{x} = \dfrac{n^n}{\left(n-1\right)!}
			\int_{x}^{\infty} t^n
			e^{-nt}\de{t}.
		\end{equation}
	Let $$\phi_n = \dfrac{\sqrt{2\pi n}\left(n/e\right)^n}{n!}.$$
	By the Stirling approximation formula, $\phi_n = 1+\Obig{n^{-1}}$ 
	for $n\to\infty$ and $\phi_n< 1$. Define 
		\begin{equation}\label{e:tau_def}
			\tau\left(t\right)=t-1-\log t \geq 0, \quad \mbox{for} 
			\quad t\in\R_+.
		\end{equation}
	In terms of $\phi_n$ and $\tau$, expression (\ref{eppsol}) acquires
the following form: 
	\begin{equation}\label{e:expp_i}
			\expp{n}{x}=\sqrt{\dfrac{n}{2\pi}} \phi_n 
			\int_{x}^{\infty}
			e^{-n\tau\left(t\right)} \de{t}.
		\end{equation}
	The integral in the right hand side can be analysed using the Laplace method. It follows from the definition that 
		\begin{equation*}
			1=\expp{n}{0} = \sqrt{\dfrac{n}{2\pi}} \phi_n 
			\int_{0}^{\infty}
			e^{-n\tau\left(t\right)} \de{t}.
		\end{equation*}
Therefore, ($\ref{e:expp_i}$) can be re-written as follows: 
		\begin{equation*}
			\expp{n}{x} = 1-\sqrt{\dfrac{n}{2\pi}} \phi_n 
			\int_{0}^{x}
			e^{-n\tau\left(t\right)} \de{t}=:1-r_n\left(x\right).
		\end{equation*}	
	Let us estimate the remainder $r_n\left(x\right)$.
	Evidently, $r_n\left(x\right)\geq 0$. An application of Taylor's theorem with 
the Lagrange form of the remainder reveals that
	for $0<t\leq x\leq 1$,
		\begin{equation}\label{e:tau_sq}
			\tau\left(t\right)\geq 
			\dfrac{\tau''\left(x\right)}{2}\left(t-x\right)^2
			+\tau'\left(x\right)\left(t-x\right)+
			\tau\left(x\right).
		\end{equation}	
	
	Noticing that $\tau'\left(x\right)=-\dfrac{1-x}{x}$ 
	and $\tau''\left(x\right)=\dfrac{1}{x^2}$ we can use the above
bound on $\tau(t)$ to obtain the following upper bound on $r_n$:
		\begin{align*}
			r_n\left(x\right)&\leq
			\sqrt{\dfrac{n}{2\pi}}\phi_n e^{-n\tau\left(x\right)} 
			\int_{0}^{x}
			e^{-\dfrac{n}{2x^2}\left(t-x\right)^2+\dfrac{n\left(1-x\right)}{x}
			\left(t-x\right)}\de{t}	
			\\
			&=
			\dfrac{\phi_n}{2}xe^{-n\left(\tau\left(x\right)-
			\left(1-x\right)^2/2\right)}
			\left(\erfc\left(\sqrt{\dfrac{n}{2}}\left(1-x\right)\right)-
			\erfc\left(\sqrt{\dfrac{n}{2}}\left(2-x\right)\right)
			\right)
			\\
			&\leq \dfrac{\phi_n}{2}xe^{-n\tau\left(x\right)}
			\erfcx\left(\sqrt{\dfrac{n}{2}}\left(1-x\right)\right),
		\end{align*}
	where $\erfc$ and $\erfcx$ are complementary and scaled 
	complementary error functions correspondingly. Finally applying 
	the classical estimate $\erfcx\left(x\right)\leq 
	\dfrac{1}{x\sqrt{\pi}}$ valid for any $x>0$ (see e.g. \cite{as}) we obtain
		\begin{equation*}
			r_n\left(x\right)\leq \dfrac{\phi_n}{\sqrt{2n\pi}} 
			\dfrac{x}{1-x}  
			e^{-n\tau\left(x\right)}< \dfrac{1}{\sqrt{2n\pi}} 
			\dfrac{x}{1-x}  
			e^{-n\tau\left(x\right)},
		\end{equation*}
	where we used that $\phi_n<1$. Therefore
		\begin{equation*}
			r_{2n}\left(1-\dfrac{\alpha_n}{2}n^{-1/2}\right)
			\leq \sqrt{\dfrac{1}{\pi}}\alpha_n^{-1} 
			e^{-2n\tau\left(1-\frac{\alpha_n}{2}n^{-1/2}\right)}.
		\end{equation*}	
	Using  \eqref{e:tau_sq} for $x=1$ and 
	$t=1-\frac{\alpha_n}{2}n^{-1/2}$ we obtain
		\begin{equation*}
			r_{2n}\left(1-\dfrac{\alpha_n}{2}n^{-1/2}\right)\leq
			\sqrt{\dfrac{1}{\pi}}\alpha_n^{-1}
			e^{-\alpha_n^2/4},
		\end{equation*}	
	which leads to the desired bound (\ref{desbnd}) for $e_{2n}^{(+)}$.
Lemma \ref{prop1} is proved.

\section{On the numerical evaluation of $p_{2n,0}$.}\label{app2}
It is clear from the proof of Lemma \ref{lm1} that the final form
of the Pfaffian or determinantal expression for the probability that an $n\times n$
real Ginibre matrix has no real eigenvalues is strongly influenced by
the choice of skew orthogonal polynomials used in the derivation. 
And even though the final exact result does not depend on the choice
of the skew orthogonal polynomials, its numerical stability is highly
sensitive to the choice.

For example, the determinantal formula (\ref{lm1f2}) is highly suitable for numerical
evaluations since the condition number of the corresponding matrix $I-M_{n}$
grows at most linearly with $n$. Indeed, its largest eigenvalue is smaller than
unity since $M_{n}$ is positive definite in virtue of the first part of Lemma \ref{lm2}. On the other hand, its smallest eigenvalue is separated
from zero by an interval of length of order $O(n^{-1})$ due to the result of
Lemma \ref{lm2} concerning the largest eigenvalue of $M_{n}$. 

This should be
contrasted to the determinantal formula derived
in \cite{KA-2005}. The condition number of the matrix
$\rho$ appearing in this formula grows exponentially with $n$, forcing one to
use high-precision numerics and leading to computation times growing
exponentially with $n$.
\end{appendix}


\begin{thebibliography}{wainrib}
 \bibitem{as} Abramowitz, M., Stegun. I. A. (eds.):
 \newblock Handbook of Mathematical Functions with Formulas, Graphs, and
 Mathematical Tables.
 \newblock {Applied Mathematics Series, vol. 55. Department of Commerce,
 National Bureau of Standards} (1964)

 \bibitem{eugene1}
 Akemann, G., Kanzieper, E.:
 \newblock Integrable structure of Ginibre’s ensemble of real
 random matrices and a Pfaffian integration theorem.
 \newblock {J. Stat. Phys.} {\bf 129}, 1159--1231 (2007)

 \bibitem{benarous}
 Ben Arous, G., Zeitouni, O.:
 \newblock Large deviations from the circular law.
 \newblock {ESAIM: Probability and Statistics} {\bf 2}, 123--134 (1998)

 \bibitem{BM:06}
 Bleher, P., Mallison, R.:
 \newblock Zeros of sections of exponential sums.
 \newblock {Int. Math. Res. Notices} {\bf 2006}, 38937 (2006)

 \bibitem{BK-2007}
 Borodin, A., Kanzieper, E.:
 \newblock A note on the Pfaffian integration theorem.
 \newblock {J. Phys. A: Math. Theor.} {\bf 40}, F849--F855 (2007)

 \bibitem{borodin_sinclair}
 Borodin, A., Sinclair, C. D.:
 \newblock The Ginibre ensemble of real random matrices and its scaling
 limits.
 \newblock {Commun. Math. Phys.} {\bf 291}, 177--224 (2009)

 \bibitem{derrida}
 Derrida, B., Zeitak, R.:
 \newblock Distribution of domain sizes in the zero temperature Glauber
 dynamics of the one-dimensional Potts model.
 \newblock {Phys. Rev. E} {\bf 54}, 2513--2525 (1996)

 \bibitem{edelman}
 Edelman, A., Kostlan, E., Shub, M.:
 \newblock How many eigenvalues of a random matrix are real?
 \newblock {J. Amer. Math. Soc.} {\bf 7}, 247--267 (1994)

 \bibitem{edelman2}
 Edelman, A.:
 \newblock The probability that a random real Gaussian matrix has $k$
 real eigenvalues, related distributions, and the circular law.
 \newblock {J. Multivariate Anal.} {\bf 60}, 203--232 (1997)

 \bibitem{forrester_nagao}
 Forrester, P. J., Nagao, T.:
 \newblock Eigenvalue statistics of the real Ginibre ensemble.
 \newblock {Phys. Rev. Lett.} {\bf 99}, 050603 (2007)

 \bibitem{forrester}
 Forrester, P. J.:
 \newblock Diffusion processes and the asymptotic bulk gap probability
 for the real Ginibre ensemble.
 \newblock arXiv:1306.4106v2 (2013)

 \bibitem{ginibre}
 Ginibre, J.:
 \newblock Statistical ensembles of complex, quaternion and real matrices.
 \newblock {J. Math. Phys.} {\bf 6}, 440--449 (1965)

 \bibitem{KA-2005}
 Kanzieper, E., Akemann, G.:
 \newblock Statistics of real rigenvalues in Ginibre’s ensemble of random
 real matrices.
 \newblock {Phys. Rev. Lett.} {\bf 95}, 230201 (2005)

 \bibitem{mcd}
 Macdonald, I. G.:
 \newblock Symmetric Functions and Hall Polynomials.
 \newblock {Clarendon Press}, Oxford (1995)

 \bibitem{ldon}
 del Molino, L. C. G., Pakdaman, K., Touboul, J., Wainrib, G.:
 \newblock The Ginibre ensemble with $k=O(n)$ real eigenvalues.
 \newblock arXiv:1501.03120v1 (2015)

 \bibitem{roger_oleg}
 Tribe, R., Zaboronski, O.:
 \newblock Pfaffian formulae for one dimensional coalescing and annihilating
   systems.
 \newblock Electron. J. Probab. {\bf 16}, 2080--2103 (2011)

 \end{thebibliography}
\end{document}